# Some applications of the Beta function

Donal F. Connon

dconnon@btopenworld.com

INTRODUCTION

The beta function $B(u,v)$ is defined for $\text{Re}\,(u) > 0$ and $\text{Re}\,(v) > 0$ by the Eulerian integral

$$B(u,v) = \int_0^1 t^{u-1}(1-t)^{v-1}\,dt$$

and it is well known that

$$B(u,v) = \frac{\Gamma(u)\Gamma(v)}{\Gamma(u+v)}$$

where $\Gamma(u)$ is the gamma function.

This work was inspired by the approach recently adopted by Morales [2].

SOME APPLICATIONS

We see that

$$B(u,v+1) = \frac{\Gamma(u)\Gamma(v+1)}{\Gamma(u+v+1)} = \frac{v\Gamma(u)\Gamma(v)}{(u+v)\Gamma(u+v)} = \frac{v}{(u+v)}B(u,v)$$

and we immediately have

$$B(u,1) = \frac{\Gamma(u)\Gamma(1)}{\Gamma(u+1)} = \frac{1}{u}$$

We have

$$\lim_{v\to 0}\frac{v}{(u+v)}B(u,v) = \lim_{v\to 0}\frac{1}{(u+v)}\lim_{v\to 0}\left[vB(u,v)\right] = \frac{1}{u}\lim_{v\to 0}\left[vB(u,v)\right]$$

but we also have $\lim_{v\to 0} B(u,v+1) = \dfrac{1}{u}$ and this implies that

$$\lim_{v\to 0}\left[vB(u,v)\right] = 1$$

We now consider

$$B(u,v)-\frac{1}{v}=\frac{vB(u,v)-1}{v}=\frac{(u+v)B(u,v+1)-1}{v}$$

and take the limit of the right-hand side

$$\lim_{v\to 0}\frac{(u+v)B(u,v+1)-1}{v}=\lim_{v\to 0}\left[(u+v)\frac{\partial}{\partial v}B(u,v+1)+B(u,v+1)\right]$$

$$=u\frac{\partial}{\partial v}B(u,v)\bigg|_{v=1}+\frac{1}{u}$$

We accordingly obtain

(1) $$\lim_{v\to 0}\left[B(u,v)-\frac{1}{v}\right]=u\frac{\partial}{\partial v}B(u,v)\bigg|_{v=1}+\frac{1}{u}$$

We have the Weierstrass canonical form of the gamma function [7, p.1]

$$\Gamma(v)=\frac{1}{v}e^{-\gamma v}\prod_{n=1}^{\infty}\left\{\left(1+\frac{v}{n}\right)^{-1}e^{\frac{v}{n}}\right\}$$

and we may write this as

$$\Gamma(v)-\frac{1}{v}=\frac{1}{v}(e^{-\gamma v}-1)\prod_{n=1}^{\infty}\left\{\left(1+\frac{v}{n}\right)^{-1}e^{\frac{v}{n}}\right\}$$

We then consider the limit

$$\lim_{v\to 0}\left[\Gamma(v)-\frac{1}{v}\right]=\lim_{v\to 0}\frac{1}{v}(e^{-\gamma v}-1)\lim_{v\to 0}\prod_{n=1}^{\infty}\left\{\left(1+\frac{v}{n}\right)^{-1}e^{\frac{v}{n}}\right\}$$

Using the Maclaurin expansion of the exponential function we readily see that

$$\lim_{v\to 0}\frac{1}{v}(e^{-\gamma v}-1)=\lim_{v\to 0}\frac{1}{v}\left[-\gamma v+O(v^2)\right]=-\gamma$$

and we then have the known result

(2) $$\lim_{v\to 0}\left[\Gamma(v)-\frac{1}{v}\right]=-\gamma$$

Alternatively we have using L'Hôpital's rule

$$\lim_{v\to 0}\left[\Gamma(v)-\frac{1}{v}\right]=\lim_{v\to 0}\left[\frac{v\Gamma(v)-1}{v}\right]=\lim_{v\to 0}\left[\frac{\Gamma(v+1)-1}{v}\right]=\lim_{v\to 0}\Gamma'(v+1)$$

and therefore we see that $\Gamma'(1)=-\gamma$.

We now consider

$$\frac{\Gamma(u+v)-\Gamma(u)}{v}=\frac{\dfrac{\Gamma(u)\Gamma(v)}{B(u,v)}-\Gamma(u)}{v}$$

$$=\Gamma(u)\frac{\Gamma(v)-B(u,v)}{vB(u,v)}$$

$$=\Gamma(u)\frac{\Gamma(v)-B(u,v)}{(u+v)B(u,v+1)}$$

$$=\Gamma(u)\frac{\Gamma(v)-\frac{1}{v}+\frac{1}{v}-B(u,v)}{(u+v)B(u,v+1)}$$

Using the definition of the derivative we have

$$\frac{d}{du}\Gamma(u)=\lim_{v\to 0}\frac{\Gamma(u+v)-\Gamma(u)}{v}$$

$$=\lim_{v\to 0}\Gamma(u)\frac{\Gamma(v)-\frac{1}{v}+\frac{1}{v}-B(u,v)}{(u+v)B(u,v+1)}$$

$$=\Gamma(u)\lim_{v\to 0}\frac{\Gamma(v)-\frac{1}{v}}{(u+v)B(u,v+1)}+\Gamma(u)\lim_{v\to 0}\frac{\frac{1}{v}-B(u,v)}{(u+v)B(u,v+1)}$$

Since $\lim_{v\to 0}(u+v)B(u,v+1)=1$ this becomes

$$\Gamma'(u) = \Gamma(u)\lim_{v\to 0}\left[\Gamma(v) - \frac{1}{v}\right] + \Gamma(u)\lim_{v\to 0}\left[\frac{1}{v} - B(u,v)\right]$$

$$= -\gamma\,\Gamma(u) - \Gamma(u)\left[u\frac{\partial}{\partial v}B(u,v)\Big|_{v=1} + \frac{1}{u}\right]$$

Hence we have an expression for the digamma function $\psi(u) = \frac{d}{du}\log\Gamma(u)$

(3) $$\frac{1}{\Gamma(u)}\frac{d}{du}\Gamma(u) = \frac{\Gamma'(u)}{\Gamma(u)} = \frac{d}{du}\log\Gamma(u) = -\gamma - \left[u\frac{\partial}{\partial v}B(u,v)\Big|_{v=1} + \frac{1}{u}\right]$$

Differentiation gives us

$$\frac{\partial}{\partial v}B(u,v) = \int_0^1 t^{u-1}\left(1-t\right)^{v-1}\log\left(1-t\right)dt$$

and we have

$$\frac{\partial}{\partial v}B(u,v)\Big|_{v=1} = \int_0^1 t^{u-1}\log\left(1-t\right)dt$$

Using (3) we obtain

(4) $$-u\int_0^1 t^{u-1}\log\left(1-t\right)dt = \gamma + \psi(u) + \frac{1}{u} = \gamma + \psi(u+1)$$

and in the case where $u$ is an integer we have the known result

$$-n\int_0^1 t^{n-1}\log\left(1-t\right)dt = \gamma + \psi(n+1) = H_n^{(1)}$$

where $H_n^{(1)}$ is the harmonic number $H_n^{(1)} = \sum_{k=1}^{n}\frac{1}{k}$.

Equation (4) may also be easily obtained from the familiar integral for the digamma function [7, p.15]

$$\psi(u+1) + \gamma = \int_0^1 \frac{1-t^u}{1-t}dt$$

where integration by parts gives us

$$= (1-t^u)\log(1-t)\Big|_0^1 - u\int_0^1 t^{u-1}\log(1-t)\,dt$$

and thus we regain

$$\psi(u+1)+\gamma = -u\int_0^1 t^{u-1}\log(1-t)\,dt$$

Using the binomial theorem we have for $|z|<1$

$$\frac{1}{(1-z)^x} = \sum_{n=0}^{\infty}\frac{(x)_n}{n!}z^n$$

where $(x)_n$ is the ascending factorial symbol (also known as the Pochhamer symbol) defined by [7, p.16] as

$$(x)_n = x(x+1)(x+2)...(x+n-1) \text{ if } n > 0 \text{ and}$$

$$(x)_0 = 1$$

It is easily seen that

$$(x)_n = \frac{\Gamma(x+n)}{\Gamma(x)}$$

We now derive an expansion for the beta function using the binomial theorem

$$B(u,v) = \int_0^1 (1-t)^{u-1}t^{v-1}dt$$

$$= \int_0^1 \sum_{n=0}^{\infty}\frac{(1-u)_n}{n!}t^n t^{v-1}dt$$

$$= \sum_{n=0}^{\infty}\frac{(1-u)_n}{(n+v)n!}$$

We therefore obtain

(5) $$B(u,v)=\frac{1}{v}+\sum_{n=1}^{\infty}\frac{(1-u)_n}{(n+v)n!}=\frac{1}{v}+\sum_{n=1}^{\infty}\frac{1}{(n+v)n!}\frac{\Gamma(n+1-u)}{\Gamma(1-u)}$$

and this gives us the limit

(6) $$\lim_{v\to 0}\left[B(u,v)-\frac{1}{v}\right]=\sum_{n=1}^{\infty}\frac{(1-u)_n}{n.n!}$$

Substituting this in (3) gives us

(7) $$\psi(u)=-\gamma-\sum_{n=1}^{\infty}\frac{(1-u)_n}{n.n!}=-\gamma-\sum_{n=1}^{\infty}\frac{1}{n.n!}\frac{\Gamma(n+1-u)}{\Gamma(1-u)}$$

$$=-\gamma-\sum_{n=1}^{\infty}\frac{1}{n^2}\frac{\Gamma(n+1-u)}{\Gamma(n)\Gamma(1-u)}=-\gamma-\sum_{n=1}^{\infty}\frac{1}{n^2B(n,1-u)}$$

with the restriction on $u$ being $1-u>0$. Therefore with $u=1/2$ we have

$$\psi\left(\frac{1}{2}\right)=-\gamma-\sum_{n=1}^{\infty}\frac{1}{n^2B(n,1/2)}$$

and from [1, p.198] we have

$$B(n,1/2)=\frac{2^{2n}}{n}\binom{2n}{n}^{-1}$$

Hence we obtain

$$\psi\left(\frac{1}{2}\right)=-\gamma-\sum_{n=1}^{\infty}\frac{1}{n2^{2n}}\binom{2n}{n}$$

and since [7, p.20] $\psi\left(\frac{1}{2}\right)=-\gamma-2\log 2$ we see that

$$\log 2=\sum_{n=1}^{\infty}\frac{1}{n2^{2n+1}}\binom{2n}{n}$$

This formula was derived in a different manner in [1, p.119].

There is a connection with the following identity

(8) $$\psi(x+a)-\psi(a)=\sum_{k=1}^{\infty}\frac{(-1)^{k+1}}{k}\frac{x(x-1)...(x-k+1)}{a(a+1)...(a+k-1)}$$

which converges for $\operatorname{Re}(x+a)>0$ . According to Raina and Ladda [5], this summation formula is due to Nörlund (see [3], [4] and also Ruben's note [6]).

Differentiation of (7) results in

(9) $$\psi'(u)=\sum_{n=1}^{\infty}\frac{1}{n^2}\frac{\Gamma(n+1-u)}{\Gamma(n)\Gamma(1-u)}[\psi(n+1-u)-\psi(1-u)]$$

and with $u=1/2$ we get

$$\psi'\left(\frac{1}{2}\right)=\sum_{n=1}^{\infty}\frac{1}{n^2}\frac{\Gamma\left(n+\frac{1}{2}\right)}{\Gamma(n)\Gamma\left(\frac{1}{2}\right)}\left[\psi\left(n+\frac{1}{2}\right)-\psi\left(\frac{1}{2}\right)\right]$$

We have Legendre's duplication formula [7, p.7] for $t>0$

$$\Gamma(t)\Gamma\left(t+\frac{1}{2}\right)=\frac{\sqrt{\pi}}{2^{2t-1}}\Gamma(2t)$$

and letting $t=n$ we obtain

$$\Gamma\left(n+\frac{1}{2}\right)=\frac{\sqrt{\pi}}{2^{2n-1}}\frac{\Gamma(2n)}{\Gamma(n)}=\frac{\sqrt{\pi}}{2^{2n-1}}\frac{(2n-1)!}{(n-1)!}=\frac{\sqrt{\pi}}{2^{2n}}\frac{(2n)!}{n!}$$

We have from [7, p.20]

$$\left[\psi\left(n+\frac{1}{2}\right)-\psi\left(\frac{1}{2}\right)\right]=2\sum_{k=1}^{n-1}\frac{1}{2k+1}$$

and we then obtain

(10) $$\psi'\left(\frac{1}{2}\right)=\sum_{n=1}^{\infty}\frac{1}{n2^{2n-1}}\frac{(2n)!}{(n!)^2}\sum_{k=1}^{n-1}\frac{1}{2k+1}$$

We have the well-known formula [7, p.22]

$$\psi^{(m)}(s)=(-1)^{m+1}m!\varsigma(m+1,s)$$

It is well known that

$$\varsigma\left(s,\frac{1}{2}\right)=\left(2^{s}-1\right)\varsigma(s)$$

and hence we obtain

$$\psi'\left(\frac{1}{2}\right)=\varsigma\left(2,\frac{1}{2}\right)=3\varsigma(2)$$

and therefore we have

(11) $$\varsigma(2)=\frac{1}{3}\sum_{n=1}^{\infty}\frac{1}{n2^{2n-1}}\frac{(2n)!}{(n!)^{2}}\sum_{k=1}^{n-1}\frac{1}{2k+1}$$

Donal F. Connon
Elmhurst
Dundle Road
Matfield
Kent TN12 7HD